\documentclass[letterpaper, 10 pt, conference]{ieeeconf}
\IEEEoverridecommandlockouts    
\usepackage[utf8]{inputenc}
\usepackage{amsmath}
\usepackage{amssymb}
\usepackage{graphicx}
\usepackage{babel}
\usepackage{balance}
\usepackage{cite}
\usepackage[position=bottom]{subfig}
\usepackage{comment}

\usepackage{mathtools}
    
\usepackage{amsmath}
\usepackage{amstext}
\usepackage{amssymb}
\usepackage{amsfonts}
\usepackage{float}

\usepackage{amsthm}  

\newtheorem{assumption}{Assumption}
\newtheorem{problem}{Problem}

\usepackage{algorithm}
\usepackage{algorithmicx} 
\usepackage{algpseudocode} %

\makeatletter
\let\NAT@parse\undefined
\makeatother
\usepackage{hyperref}
\hypersetup{
   colorlinks=true,
    linkcolor= blue,
    allcolors=blue,
    citecolor = blue,
    filecolor=black,      
    urlcolor=blue,
    }

\title{\LARGE \bf Congestion-Aware Routing, Rebalancing, and Charging Scheduling for Electric Autonomous Mobility-on-Demand System}

\author{Heeseung Bang, \textit{Student Member, IEEE}, Andreas A. Malikopoulos, \textit{Senior Member, IEEE}
    \thanks{This research was supported by the Sociotechnical Systems Center (SSC) at the University of Delaware.}
	\thanks{The authors are with the Department of Mechanical Engineering, University of Delaware, Newark, DE 19716, USA. (emails: \tt\small{heeseung@udel.edu}; \tt\small{andreas@udel.edu}.)}
}

\date{August 2021}

\begin{document}

\maketitle

\begin{abstract}
In this paper, we investigate the problem of routing, rebalancing, and charging for electric autonomous mobility-on-demand systems concerning traffic congestion. We analyze the problem at the macroscopical level and use a volume-delay function to capture traffic congestion. To address this problem, we first formulate an optimization problem for routing and rebalancing. Then, we present heuristic algorithms to find the loop of the traffic flow and examine the energy constraints within the resulting loop. We impose charging constraints on the re-routing problem so that the new solution satisfies the energy constraint.
Finally, we verify the effectiveness of our method through simulation.
\end{abstract}

\section{Introduction}

\PARstart{O}{ver} the past few decades, the population of cities around the world has been rapidly growing. Currently, more than 55$\%$ of the world's population lives in cities, and it is predicted to increase to 68$\%$ by 2050 \cite{urbanization2018}. According to the Texas A\&M Transportation Institute \cite{Schrank2019}, traffic congestion has dramatically increased over the last decade. In 2017, an average commuter in the US spent an extra $54$ hours and $21$ gallons of fuel due to traffic congestion, corresponding to \$$179$ billion. 
Autonomous vehicles mark a paradigm shift in which a myriad of opportunities exist for users to better monitor the transportation network conditions and make optimal operating \cite{zhao2019enhanced} and logistics decisions \cite{remer2019multi}.


One of the promising ways to alleviate traffic issues is to utilize an autonomous mobility-on-demand (AMoD) system.
AMoD system can efficiently manage vehicles to provide a ride to the customers with fewer vehicles, which reduces vehicles on the roads.
In addition, since one can centrally control AMoD system, it is possible to consider a system-optimal decision.
In this paper, we consider an electric autonomous mobility-on-demand (E-AMoD) system, in which battery constraints should also be considered.

There are two areas related to an E-AMoD system: (1) routing for AMoD and (2) charging scheduling for electric vehicles. 
The routing problem considering traffic congestion started from the work of Salazar et al. \cite{salazar2019congestion}, where they introduced a congestion-aware routing scheme with a piecewise-affine approximation of the travel delay function to solve the problem in real-time.
Wollenstein-Betech et al. \cite{wollenstein2020congestion} considered the endogenous impact of the private vehicles' routing decision by formulating a game-theoretical problem.
They extended the work \cite{salazar2019congestion} and \cite{wollenstein2020congestion}, and presented a route-recovery algorithm to find the exact route for specific travel demands in \cite{wollenstein2021routing}.
There has also been a series of papers considering microscopic-level routing for different objectives, e.g., ride-sharing \cite{tsao2019model} and vehicle-trip assignment \cite{chen2020optimal}.
Smith et al. \cite{smith2013rebalancing} introduced a rebalancing strategy for mobility-on-demand system vehicles and drivers by using a fluid-approximated model. 
A more comprehensive review of the research efforts on AMoD is provided in \cite{zardini2021analysis}. There are also several review papers providing a good summary of shared mobility, e.g., see \cite{Jorge2013, Agatz2012a, Furuhata2013a, Brandstatter2016, lavieri2017, jittrapirom2017, utriainen2018}.

Various research efforts have been addressed in the literature, which considered routing and charging scheduling problems for E-AMoD system by using different approaches, e.g., deep reinforcement learning \cite{liang2020mobility,wan2018model}, mixed-integer linear programming \cite{boewing2020vehicle,chen2019integrated,salazar2019optimal}, and a graph search algorithm \cite{sachenbacher2011efficient}.
Rossi et al. \cite{rossi2019interaction} developed a network flow and considered fleet impact on the power network.
Estandia et al. \cite{estandia2021interaction} extended this work and additionally considered the power distribution network.
Li et al. \cite{li2019agent} presented a vehicle relocation strategy for E-AMoD by matching nearest idle vehicles. If no vehicle can pick up the customer in time, the customer leaves the system.
Pourazarm et al. \cite{pourazarm2014optimal} presented an optimal routing method for electric vehicles considering battery constraints as imposed by the vehicles components' sizing \cite{Malikopoulos2013b} and consumers' preferences \cite{Shaltout2014}. 
They used a dynamic programming approach, which was effective only for single-vehicle routing.

In this paper, we consider the congestion-aware routing, rebalancing, and charging scheduling for E-AMoD. We formulate the optimization problem to find the best route for traveling, rebalancing, and charging. The objective function consists of the customers' travel time and rebalancing cost.
Although there have been a lot of efforts on charging scheduling for electric vehicles as well as routing for AMoD system, to the best of our knowledge, considering both in conjunction with traffic congestion has not yet been addressed in the literature. 
The main contributions of this paper are: (1) an approach of evaluating the energy constraint of E-AMoD vehicles at the macroscopic level; and (2) an approach to combine charging scheduling to the congestion-aware routing problem.

The remainder of the paper is organized as follows. In Section \ref{sec:problem}, we present our model and problem formulation. In Section \ref{sec:solution}, we solve the routing and rebalancing problem and then introduce the approaches for energy constraint evaluation and re-planning. We verify our method through simulations in Section \ref{sec:simulation}, and provide concluding remarks in Section \ref{sec:conclusion}.

\section{Problem Formulation} \label{sec:problem}

We consider a macroscopic model for routing, rebalancing, and charging scheduling. Rather than focusing on each vehicle, we formulate the problem using traffic flow.
The road network is denoted by a directed graph $\mathcal{G} = (\mathcal{V},\mathcal{E})$, where $\mathcal{V}\subset\mathbb{N}$ is a set of intersection (vertices) and $\mathcal{E}\subset\mathcal{V}\times\mathcal{V}$ is a set of streets (edges) in the road network.
Let $\mathcal{C}=\{c_1,\dots,c_N\}\subset\mathcal{V}$ be a set of $N\in\mathbb{N}$ charging stations.

We consider $M\in\mathbb{N}$ travel demands. For each travel demand $ m\in\mathcal{M} = \{1,\dots,M\}$, let $o_m \in\mathcal{O}$, $d_m\in\mathcal{D}$, and $\alpha_m\in\mathbb{R}_{>0}$ denote origin, destination, and demand rate, respectively, where $\mathcal{O}\subset\mathcal{V}$ is a set of origins, and $\mathcal{D}\subset\mathcal{V}$ is a of destinations. Demand rate is defined as the number of customers per unit time.
For each road $(i,j)\in\mathcal{E}$, let $u_{ij}^m$ denote \textit{customer’s flow} traveling from $o_m$ to $d_m$, $r_{ij}^{ab}$ denote \textit{rebalancing flow} that relocates the vehicles from $a\in\mathcal{D}$ to $b\in\mathcal{O}$, and $p_{ij}$ denote \textit{private vehicle’s flow}.
One can consider private vehicles as self-interested game players and assume some rationale behind their decisions as reported in \cite{wollenstein2020congestion}. However, this is beyond the scope of this paper; hence we consider private vehicle’s flow to be constant.
Let
\begin{align}
    u_{ij} &= \sum_{m\in\mathcal{M}} u_{ij}^m,~~\forall(i,j)\in\mathcal{E},\\
    r_{ij} &= \sum_{a\in\mathcal{D}}\sum_{b\in\mathcal{O}} r_{ij}^{ab},~~\forall(i,j)\in\mathcal{E},
\end{align}
where $u_{ij}$ is the total customer's flow and $r_{ij}$ is total rebalancing flow on the road $(i,j)\in\mathcal{E}$.
Then, the total flow on each road becomes $x_{ij} = u_{ij} + r_{ij} + p_{ij},~~\forall(i,j)\in\mathcal{E}$.
To meet all demand, the customer's flow must satisfy
\begin{align}
    &\sum_{i:(i,j)\in\mathcal{E}}u_{ij}^m = \sum_{k:(j,k)\in\mathcal{E}}u_{jk}^m,~\forall m\in\mathcal{M},j\in\mathcal{V}\setminus\{o_m,d_m\},\label{eqn:con_u1}\\
    & \sum_{k:(j,k)\in\mathcal{E}}u_{jk}^m = \alpha_m,~\forall m\in\mathcal{M},j=o_m,\label{eqn:con_u2}\\
    & \sum_{i:(i,j)\in\mathcal{E}}u_{ij}^m = \alpha_m,~\forall m\in\mathcal{M},j=d_m.\label{eqn:con_u3}
\end{align}
Constraint \eqref{eqn:con_u1} takes care of flow conservation on each node, while constraints \eqref{eqn:con_u2} and \eqref{eqn:con_u3} balance the customer's flow with the demand rate $\alpha_m$ and guarantee the flow starts and ends at $o_m$ and $d_m$, respectively.
Next, we define the net rebalancing flow.
For each node $j\in\mathcal{O}\cup\mathcal{D}$, the {\it net rebalancing flow} is given by $r_j^{\text{net}} = \sum_{m\in\mathcal{M}} \left(\mathbf{1}_{j=d_m}- \mathbf{1}_{j=o_m} \right) \alpha_m$, where $\mathbf{1}_{j=h}$ is an indicator function whose value is $1$ if $j=h$ and $0$ otherwise.

The net rebalancing flow describes, at node $j\in\mathcal{O}\cup\mathcal{D}$, how much rebalancing flow needs to go out if $r_j^{\text{net}}>0$, or come in if $r_j^{\text{net}}<0$.
To ensure the rebalancing flow meets the net flow, we impose the following constraints 
\begin{align}
    \sum_{a\in\mathcal{D}}\sum_{i:(i,b)\in\mathcal{E}} r_{ib}^{ab} &= \begin{cases} 0, ~&~r_b^{\text{net}} > 0\\ -r_b^{\text{net}}, ~&~r_b^{\text{net}} \leq 0\end{cases},~\forall b\in\mathcal{O}, \label{eqn:con_u4}\\
    \sum_{b\in\mathcal{O}}\sum_{k:(a,k)\in\mathcal{E}} r_{ak}^{ab} &= \begin{cases} r_a^{\text{net}}, ~&~r_a^{\text{net}} > 0\\ 0, ~&~r_a^{\text{net}} \leq 0\end{cases},~\forall a\in\mathcal{D}, \label{eqn:con_u5}\\
    \sum_{a\in\mathcal{D}}\sum_{k:(b,k)\in\mathcal{E}} r_{bk}^{ab} &= 0,~\forall b\in\mathcal{O}, \label{eqn:con_u6}\\
    \sum_{b\in\mathcal{O}}\sum_{i:(i,a)\in\mathcal{E}} r_{ia}^{ab} &= 0,~\forall a\in\mathcal{D}. \label{eqn:con_u7}
\end{align}
The constraints \eqref{eqn:con_u4} -- \eqref{eqn:con_u7} ensure all  nodes to have net rebalancing flow without any redundant flow. 
Lastly, we impose flow conservation for rebalancing flow, which can be expressed by
\begin{equation}
    \sum_{i:(i,j)\in\mathcal{E}} r_{ij}^{ab}=\sum_{k:(j,k)\in\mathcal{E}} r_{jk}^{ab},~\forall a\in\mathcal{D},b
    \in\mathcal{O},j\in\mathcal{V}\setminus\{a,b\}. \label{eqn:con_r}
\end{equation}

We consider congestion on the road, by defining the travel time function $t_{ij}(x)$ for each road $(i,j)\in\mathcal{E}$.
Generally, the travel time function has the following structure \cite{spiess1990technical}, i.e., $t_{ij}(x_{ij}) = t_{ij}^0 f\left(\frac{x_{ij}}{\gamma_{ij}}\right)$.
It is a product of free-flow travel time $t_{ij}^0\in\mathbb{R}_{>0}$ and a normalized congestion function, where $\gamma_{ij}$ is the capacity of the road. The most widely used travel time function is the one proposed by \textit{Bureau of Public Roads} (BPR) \cite{us1964traffic}, which is also called BPR function.
By using the BPR function, the travel time can be expressed as
\begin{equation}
    t_{ij}(x_{ij}) = t_{ij}^0\left(1+0.15\left(\frac{x_{ij}}{\gamma_{ij}}\right)^4\right). \label{eqn:BPR}
\end{equation}

Recall that our goal is to find the optimal route and flow for traveling, rebalancing, and charging scheduling for E-AMoD system.
The objective is to minimize the total travel time of the customer's flow and the rebalancing cost.
Therefore, we formulate the routing and rebalancing problem as follows.

\begin{problem} \label{prb:route}
We find the routes for traveling and rebalancing flows by solving the following optimization problem:
\begin{equation}
\begin{aligned}
    \min_{\mathbf{u},\mathbf{r}} ~&J(\mathbf{u},\mathbf{r}) = \sum_{(i,j)\in\mathcal{E}} \bigg\{ t_{ij}(x_{ij})u_{ij} + w_r r_{ij}\bigg\} \label{eqn:optimization}\\
    \text{s.t. } & \eqref{eqn:con_u1}\text{ -- }\eqref{eqn:con_u3}, \eqref{eqn:con_u4}\text{ -- }\eqref{eqn:con_r},\\
    &\mathbf{u} \geq 0, \mathbf{r} \geq 0,
\end{aligned}    
\end{equation}
where $w_r\in\mathbb{R}_{>0}$ is weighted cost for rebalancing, $\mathbf{u}$ is a vector of $u_{ij}^m$ for all $m\in\mathcal{M}, (i,j)\in\mathcal{E}$, and $\mathbf{r}$ is a vector of $r_{ij}^{ab}$ for all $a\in\mathcal{D}, b\in\mathcal{O}, (i,j)\in\mathcal{E}$.
\end{problem}


In our modeling framework, we impose the following assumptions.

\begin{assumption}\label{smp:macroscopic}
We neglect microscopic phenomena such as intersection delay due to signal lights and traffic delay caused by picking up and dropping off customers.
\end{assumption}
\begin{assumption}\label{smp:timeinvariance}
We assume that the travel demands are time-invariant.
\end{assumption}
\begin{assumption}\label{smp:chargingstation}
There are enough charging stations such that any vehicle at any location can visit charging stations with 10\% of the battery.
\end{assumption}
\begin{assumption}\label{smp:chargingspace}
We assume that there are always enough chargers at the charging station so that the balancing between incoming flow and outgoing flow of the charging station remains the same.
\end{assumption}

Assumptions \ref{smp:macroscopic} and \ref{smp:timeinvariance} are acceptable for macroscopic analysis, especially in an urban area where the heavy amount of traveling seriously affects traffic congestion.
Assumption \ref{smp:chargingstation} may be strong, but we can relax this assumption by considering the fleet size or finding the best state of charge (SoC) to reroute to a charging station.
Assumption \ref{smp:chargingspace} may be relaxed by adding a flow limit on the charging station.

The energy consumption of electric vehicles depends on the speed. The optimal speed for energy-efficient driving is in the range of 10 to 50 km/h \cite{gree2020cloud}. By using the relationship between energy usage and speed of the vehicle, we compute the energy required for traveling the road $(i,j)\in\mathcal{E}$ as
\begin{equation}
    E_{ij} = s_{ij} E_f\left(\frac{s_{ij}}{t_{ij}(x_{ij})}\right), \label{eqn:energy}
\end{equation}
where $s_{ij}$ is the length of the road $(i,j)\in\mathcal{E}$ and $E_f$ is an energy function that maps the average speed of the vehicle to the energy consumption. We compute the average speed on the road using the length $s_{ij}$ and the travel time function \eqref{eqn:BPR}.
As shown in \eqref{eqn:energy}, the energy required to travel the road $(i,j)\in\mathcal{E}$ is a function of the flow $x_{ij}$. Hence, once we get the solution to the Problem \ref{prb:route}, we can directly get the required energy for all the roads.

\section{Solution Approach} \label{sec:solution}
In this section, we (1) present our solution approach for the optimization problem, (2) explain how we examine whether or not the current route meets the energy constraints, and (3) formulate a re-routing problem when energy constraints are not satisfied.

\subsection{Convex Approximation}
The problem of routing and rebalancing is non-convex, because of the product term $t(x_{ij})u_{ij}$ in the objective function. This non-convexity significantly increases the computational cost of the problem. Therefore, we find the upper bound of the objective function, which is convex, so that we can solve the convex optimization problem. 

\begin{figure}[t]
    \centering
    \includegraphics[scale=0.13]{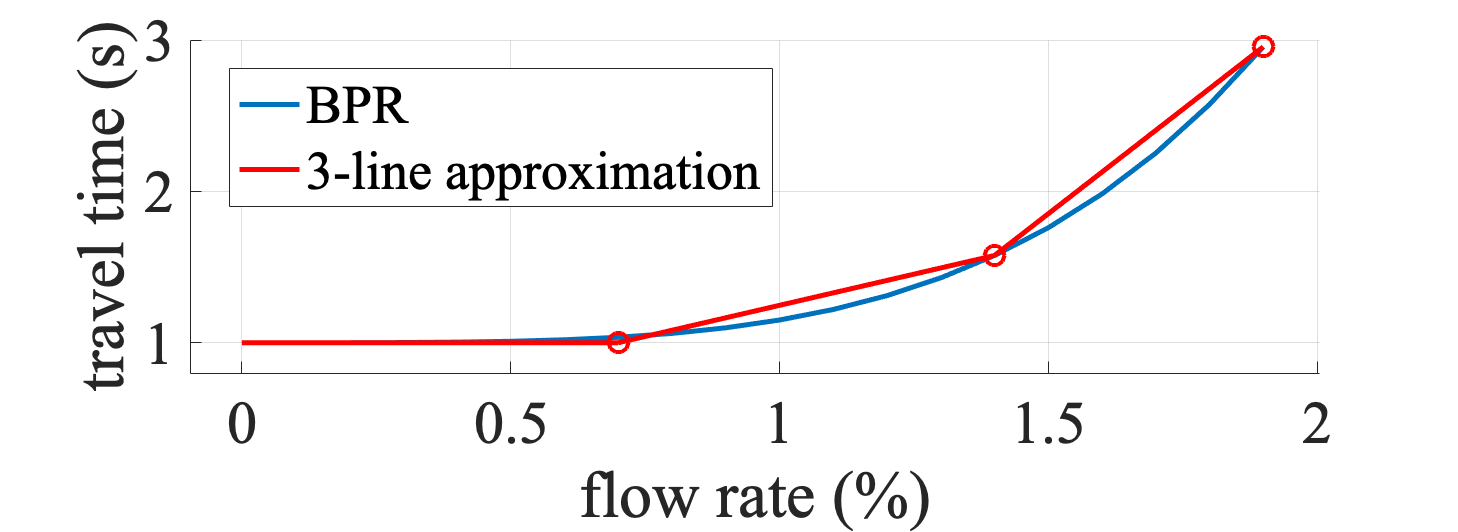}
    \caption{BPR function (blue) and 3-line affine approximation (red).}
    \label{fig:BPR}
    
\end{figure}[t]

We use the same approach as presented in \cite{salazar2019congestion,wollenstein2020congestion,wollenstein2021routing}. In order to find the upper bound of the objective function, we first use a 3-line piecewise affine function to approximate the travel time function as shown in Fig. \ref{fig:BPR}.
Let $\Tilde{t}_{ij}(x_{ij})$ be the following piecewise-affine function
\begin{align}
    &\Tilde{t}_{ij}(x_{ij}) \nonumber\\
    &=\begin{cases}
    t_{ij}^0, ~~~~~~~~~~~~~~~~~~~~~~~~~~~~~~~~~~~~~~~~~~\text{ if } x_{ij} < x^{\text{th}1}_{ij}, \\
    t_{ij}^0\left(1+a_{ij}\left(x_{ij}-x^{\text{th}1}_{ij}\right)\right), ~~~~~~~\text{ if } x^{\text{th}1}_{ij} \leq x_{ij} < x^{\text{th}2}_{ij}, \\
    t_{ij}^0\left(1+a_{ij}\left(x^{\text{th}2}_{ij}-x^{\text{th}1}_{ij}\right)+b_{ij}\left(x_{ij}-x^{\text{th}2}_{ij}\right)\right),\\
    ~~~~~~~~~~~~~~~~~~~~~~~~~~~~~~~~~~~~~~~~~~~~~\text{ if } x^{\text{th}2}_{ij} \leq x_{ij},
    \end{cases}
\end{align}
where $x^{\text{th}1}_{ij}$ and $x^{\text{th}2}_{ij}$ are the end points of the first two line segments, respectively, and $a_{ij} = l_1/\gamma_{ij}$, $b_{ij}=l_2/\gamma_{ij}$. Here, $l_1$ and $l_2$ are the slope of the second and third segments in Fig. \ref{fig:BPR}, respectively.
Next, in order to simplify the piecewise-affine function, we define the slack variables as
$\varepsilon_{ij}^{(1)} = \max\{0,x_{ij}-x^{\text{th}1}_{ij}-\varepsilon_{ij}^{(2)}\}, 
\varepsilon_{ij}^{(2)} = \max\{0,x_{ij}-x^{\text{th}2}_{ij}\}$.
The slack variable $\varepsilon_{ij}^{(1)}\in[0,x^{\text{th}2}_{ij}-x^{\text{th}1}_{ij}]$ is the excess flow over $x^{\text{th}1}_{ij}$, and $\varepsilon_{ij}^{(2)}$ is the excess flow over $x^{\text{th}2}_{ij}$.
By using these slack variables, the piecewise-affine function can be expressed as 
$\Tilde{t}_{ij}(x_{ij}) = t_{ij}^0 (1 + a_{ij}\varepsilon_{ij}^{(1)} + b_{ij}\varepsilon_{ij}^{(2)})$.
Thus, we can define an element-wise approximation $\Tilde{J}_{ij}$ of the objective function, i.e.,

\begin{subequations}
\begin{align}
    \Tilde{J}_{ij} &= \Tilde{t}_{ij}(x_{ij})u_{ij} + w_r r_{ij}\\
    &= t_{ij}^0 (1 + a_{ij}\varepsilon_{ij}^{(1)} + b_{ij}\varepsilon_{ij}^{(2)})u_{ij} + w_r r_{ij}\label{eqn:sub1}\\
    &~\begin{aligned} =~ &t_{ij}^0u_{ij} + t_{ij}^0 a_{ij}\varepsilon_{ij}^{(1)}(\varepsilon_{ij}^{(1)}+\varepsilon_{ij}^{(2)}+x_{ij}^{\text{th}1}-r_{ij}-p_{ij})\\
    &+t_{ij}^0 b_{ij}  \varepsilon_{ij}^{(2)}(\varepsilon_{ij}^{(2)}+x_{ij}^{\text{th}2}-r_{ij}-p_{ij})+ w_r r_{ij}.
    \end{aligned}\label{eqn:sub2}
\end{align}
\end{subequations}

Next, we find the upper bound of $\Tilde{J}_{ij}$ and minimize it to find an approximated solution. To find the convex upper bound, we remove all the non-convex term, which are product of two optimization variables such as $\varepsilon_{ij}^{(1)}\varepsilon_{ij}^{(2)}$, $\varepsilon_{ij}^{(1)}r_{ij}$, and $\varepsilon_{ij}^{(2)}r_{ij}$. We can simply change $\varepsilon_{ij}^{(1)}\varepsilon_{ij}^{(2)}$ into $(x^{\text{th}2}_{ij}-x^{\text{th}1}_{ij})\varepsilon_{ij}^{(2)}$ , since $\varepsilon_{ij}^{(2)} = 0$ if $x_{ij} < x_{ij}^{\text{th}2}$ and $\varepsilon_{ij}^{(1)} = (x_{ij}^{\text{th}2}-x_{ij}^{\text{th}1})$ if $x_{ij}>x_{ij}^{\text{th}2}$ and 0 otherwise.
By adding $(a_{ij}t_{ij}^0\varepsilon_{ij}^{(1)}+b_{ij} t_{ij}^0 \varepsilon_{ij}^{(2)})r_{ij}$, we attain the upper bound of $\Tilde{J}_{ij}$, which can be written as following quadratic function

\begin{align}
    J^{QP}_{ij} =~ & t_{ij}^0 u_{ij}+a_{ij}t_{ij}^0\varepsilon_{ij}^{(1)}(\varepsilon_{ij}^{(1)}+x^{\text{th}1}_{ij}-p_{ij})\\
    &+ b_{ij}t_{ij}^0\varepsilon_{ij}^{(2)}(\varepsilon_{ij}^{(2)}+x^{\text{th}2}_{ij}-p_{ij})\nonumber\\
    &+a_{ij}t_{ij}^0 (x^{\text{th}2}_{ij}-x^{\text{th}1}_{ij})\varepsilon_{ij}^{(2)} + w_r r_{ij}.\nonumber
\end{align}
From the definition, slack variables have following inequality constraints
\begin{align}
    \varepsilon_{ij}^{(1)} \geq 0,~&~\varepsilon_{ij}^{(1)}\geq x_{ij}-x^{\text{th}1}_{ij}-\varepsilon_{ij}^{(2)}, \label{eqn:slack1}\\
    \varepsilon_{ij}^{(2)} \geq 0,~&~\varepsilon_{ij}^{(2)}\geq x_{ij}-x^{\text{th}2}_{ij}. \label{eqn:slack2}
\end{align}
The approximation of solution to the Problem \ref{prb:route} can be obtained by solving the problem defined next.

\begin{problem} \label{prb:convex}
    To find an approximated solution, we solve quadratic program defined as
    \begin{align}
    \min_{\textbf{u},\textbf{r}} ~&~ \sum_{(i,j)\in\mathcal{E}} J_{ij}^{QP}\\
    \text{s.t. } & \eqref{eqn:con_u1}\text{ -- }\eqref{eqn:con_u3}, \eqref{eqn:con_u4}\text{ -- }\eqref{eqn:con_r}, \eqref{eqn:slack1},\eqref{eqn:slack2}. \nonumber
\end{align}
\end{problem}

\subsection{Energy Consumption of Flow}
In order to examine if the vehicles need to be charged or not, we require the knowledge of the vehicles' next trip. Therefore we form loops for the vehicles and find the traffic flow of each loop by running Algorithm \ref{alg:loop}.

\begin{algorithm}
\caption{Loop Recovery}\label{alg:loop}
\begin{algorithmic}[1]

\While{there is unlooped flow}
    \State{Select any $v\in\mathcal{O}\cup\mathcal{D}$ that has unlooped flow}
    \While{True}
        \State{Check $v$ as visited node}
        \State{Find all the nodes connected to $v$}
        \If{there is node that is already visited}
            \State{Generate loop}   
            \State{Break}
        \Else
            \State{$v \gets$ connected node with maximum flow}
        \EndIf
    \EndWhile
    \State{Remove the loop from unlooped flow}
\EndWhile
\end{algorithmic}
\end{algorithm}

Next, we present Algorithm \ref{alg:chargingschedule}, which explains the method of charging scheduling. Since E-AMoD vehicles recurrently travel along the loop, each loop requires at least one charging schedule.
We schedule it at the beginning of the loop to ensure that the vehicles do not run out of battery if they can safely travel a single loop.
In addition, if there is any rebalancing trip in the loop, we select it to be the first trip and have the first charging schedule because the charging trip between the customers' travel requires additional flow on the roads.
For the rest of the trip in the loop, we compute energy consumption and estimate SoC of the vehicles by using Algorithms \ref{alg:energy} and \ref{alg:soc}, which will be explained later in this subsection.
Because of Assumption \ref{smp:chargingstation}, vehicles are safe to travel until the SoC subceeds $10\%$ of the total battery capacity.
If the SoC is expected to be lower than $10\%$, we schedule charging during the rebalancing trip.

\begin{algorithm}
\caption{Charging Scheduling}\label{alg:chargingschedule}
\begin{algorithmic}[1]
\For{each loop}
    \State{$x_{CS} \gets $ flow of the loop}
    \If{there is rebalancing flow in the loop}
        \State{Select rebalancing as first trip}
    \Else
        \State{Select any customer's travel as first trip}
    \EndIf
    \For{each trip}
        \State{$(o,d) \gets $ origin and destination of the trip}
        \State{$E \gets $ ENERGY(trip)}
        \If{$E_\mathrm{bat} - E < 10\%$}
            \State{$E_\mathrm{bat} \gets$ SoC(trip, $E$)}
            \If{trip is rebalancing}
                \State{$\mathcal{S} \gets \mathcal{S}\cup\{(o,d,x_{CS})\}$}
            \Else
                \State{$\mathcal{S} \gets \mathcal{S}\cup\{(o,o,x_{CS})\}$}
            \EndIf
        \Else
            \State{$E_\mathrm{bat} \gets E_\mathrm{bat} - E$}
        \EndIf
    \EndFor
    \If{there exists same origin and destination $(o,d)$}
        \State{$x_{\mathrm{sum}} \gets$ accumulate all the flow for $(o,d)$}
        \State{Remove all tuples that has $(o,d)$ from $\mathcal{S}$}
        \State{$\mathcal{S} \gets \mathcal{S}\cup(o,d,x_{\mathrm{sum}})$}
    \EndIf
\EndFor

\end{algorithmic}
\end{algorithm}

\begin{algorithm}
\caption{Computation of Energy Consumption}\label{alg:energy}
\begin{flushleft}
        \textbf{Input:} trip, $(o,d)$  \\
        \textbf{Output:} $E$
\end{flushleft}
\begin{algorithmic}[1]
\State{$E_{c}(o) \gets 0$;~$i \gets o$}
\While{$i \neq d$}
    \State{Check $i$ as visited}
    \For{each node $j$ connected to $i$ in the trip}
        \State{Compute $E_{ij}$ from \eqref{eqn:energy}}
        \State{$E_c(j) \gets \max\{E_c(j), E_c(i) + E_{ij}\}$}
        \State{Disconnect the flow from $i$ to $j$ in the trip}
    \EndFor
    \State{$i \gets$ unvisited node with no incoming flow}
\EndWhile
\State{$E \gets E_{c}(d)$}
\end{algorithmic}
\end{algorithm}

Algorithm \ref{alg:energy} computes the maximum required energy for the given trip.
Note that there could be multiple routes for each trip.
To avoid the situation of running out of battery, we compute the maximum energy among the different routes.
The algorithm is similar to the breadth-first search in a directed graph, but the iteration stops only after searching all the routes.
At every node, we compute the required energy for traveling to all the connected nodes by using \eqref{eqn:energy}.
Whenever we find the higher required energy, we select it to be the required energy for the trip.

\begin{algorithm}
\caption{Estimation of State of Charge}\label{alg:soc}
\begin{flushleft}
        \textbf{Input:} trip, $E$, $(o,d)$  \\
        \textbf{Output:} $E_\mathrm{bat}$
\end{flushleft}
\begin{algorithmic}[1]
\If{$d\in\mathcal{C}$}
    \State{$E_\mathrm{bat} \gets 100\%$}
\ElsIf{$o\in\mathcal{C}$}
    \State{$E_\mathrm{bat} \gets 100\%-E$}
\Else
    \If{trip is rebalancing}
        \State{$E_\mathrm{bat} \gets 90\%$}
    \Else
        \State{$E_\mathrm{bat} \gets 90\% - E$}
    \EndIf
\EndIf
\end{algorithmic}
\end{algorithm}

Algorithm \ref{alg:soc} presents the estimated SoC after charging.
We first consider the rebalancing trip.
If the destination of the trip has charging station, SoC becomes $100\%$ as vehicles can charge at the destination.
Meanwhile, if the origin of the trip has charging station, vehicles consume the energy required for the trip after charging; hence, it becomes $100\% - E$. In fact, one can schedule charging during the rebalancing trip to keep SoC high at the destination. However, the difference in SoC may be trivial, while the variation of the flow due to the charging schedule may increase the cost. Therefore, we avoid changing the flow if possible.
If neither the origin nor destination has a charging station, vehicles need to visit a charging station during the rebalancing trip.
We estimate SoC to be $90\%$ as a trip to the charging station could require up to $10\%$.
If it is the customer's trip, the vehicle cannot visit a charging station during the trip.
Thus, they need to schedule extra trip, which departs from the origin, visits a charging station, and returning to the origin of the trip.

\subsection{Re-Routing with Consideration of Charging}
After scheduling charging through Algorithm \ref{alg:chargingschedule}, we impose additional constraints to ensure that the flow visits charging stations.

\begin{align}
    \sum_{k:(a,k)\in\mathcal{E}} r_{ak}^{aa} & = x_{CS},~ \forall (a,a,x_{CS}) \in\mathcal{S},  \label{eqn:charge1}\\
    \sum_{i:(i,a)\in\mathcal{E}} r_{ia}^{aa} & = x_{CS},~ \forall (a,a,x_{CS}) \in\mathcal{S}, \label{eqn:charge2}\\
    \sum_{c\in\mathcal{C}}\sum_{i:(i,c)\in\mathcal{E}} r_{ic}^{aa} & \geq x_{CS},~ \forall (a,a,x_{CS})\in\mathcal{S}, \label{eqn:charge3} \\
    \sum_{c\in\mathcal{C}}\sum_{i:(i,c)\in\mathcal{E}} r_{ic}^{ab} & \geq x_{CS},~\forall (a,b,x_{CS})\in\mathcal{S}. \label{eqn:charge4}
\end{align}
Here, $\mathcal{S}$ is a set of charging schedules, which is obtained from Algorithm \ref{alg:energy}.
It consists of tuples with origin-destination and the charging flow $x_{CS}$.
In case the charging is scheduled before the customer's travel, \eqref{eqn:charge1} -- \eqref{eqn:charge3} ensures the flow visits the charging station before the vehicles pick up the customers.
On the other hand, if the charging is scheduled in the rebalancing trip, constraint \eqref{eqn:charge4} simply lets the flow visit the charging station during rebalancing.
By using these charging constraints, we define the re-routing problem.

\begin{problem} \label{prb:reroute}
For re-routing problem we solve the following optimization:
\begin{align}
    \min_{\textbf{u},\textbf{r}} ~&~ \sum_{(i,j)\in\mathcal{E}} J_{ij}^{QP}\\
    \text{s.t. } & ~\eqref{eqn:con_u1} \text{ -- } \eqref{eqn:con_r}, \eqref{eqn:slack1},\eqref{eqn:slack2},\eqref{eqn:charge1} \text{ -- }\eqref{eqn:charge4}. \nonumber
\end{align}
\end{problem}


\section{Simulation Results} \label{sec:simulation}

In this section, we present the simulation results to evaluate our approach and analyze how much charging trip affects the flow.
We considered a small road network with 8 nodes and 22 edges.
We randomly selected private vehicles' flow and fixed it for different cases.



\begin{figure*}[h]
    \centering
    \includegraphics[width=\linewidth]{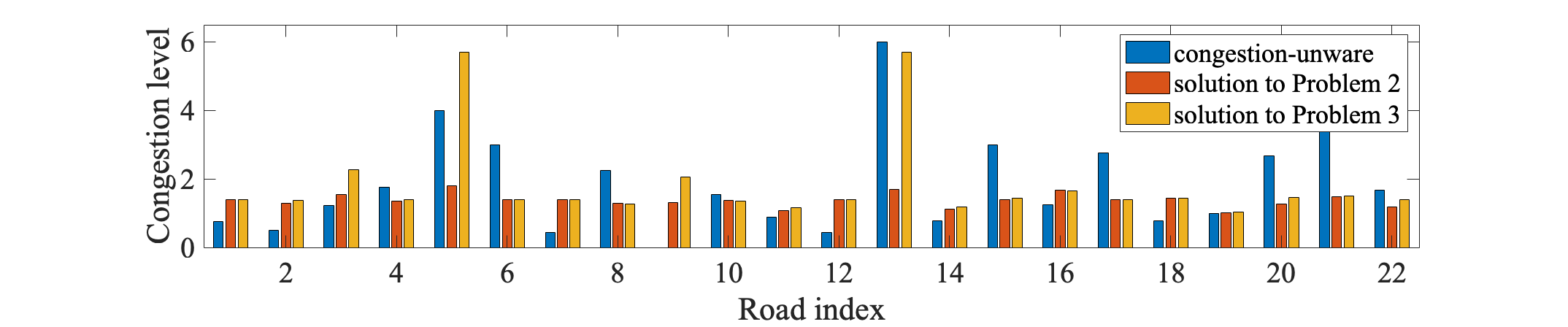}
    \caption{Comparison of congestion level on each road for difference approaches.}
    \label{fig:congestion}
\end{figure*}

To demonstrate the effectiveness of our approach, we conducted a simulation for the congestion-unaware approach and compared the result. We used free-flow travel time for each road and solved Problem \ref{prb:convex}.
Figure \ref{fig:congestion} shows the difference in congestion level for all the solution approaches on each road.
The congestion-unaware approach relies on the roads known to have short free-flow travel time even if the roads have a high level of traffic congestion.
Meanwhile, the method barely used certain roads, e.g., road 9 in Fig. \ref{fig:congestion} had no flow at all. On the other hand, the congestion-aware solution obtained from Problem \ref{prb:convex} distributed all the flow and used all the roads, which resulted in remaining the flow less than twofold of the road capacity.
Considering charging scheduling, additional trips for charging dramatically increased the flow on the roads 5 and 13 in Fig. \ref{fig:congestion}.
However, the traffic congestion caused by additional trips may be inevitable since charging is mandatory for vehicles.





\section{Concluding Remarks} \label{sec:conclusion}
This paper proposed a solution approach to congestion-aware routing problems considering traveling, rebalancing, and charging scheduling for E-AMoD system. The idea is to solve the routing problem,  examine energy constraints along the given path, and re-route in case charging is needed.
We validated the effectiveness of our approach through simulation.

Ongoing work considers microscopic phenomena using connected and automated vehicles \cite{chalaki2020TCST,chalaki2020TITS,Malikopoulos2020} in a scaled environment \cite{Beaver2020DemonstrationCity}. Future work should consider incorporating the energy constraints into the optimization problem, as well as finding an optimal loop in terms of charging scheduling. 
\balance 

\bibliographystyle{IEEEtran}
\bibliography{Bang, IDS, SharedMobilityRef}

\end{document}